\documentclass[12pt]{article}
\usepackage{amssymb,amsmath}
\usepackage{cases}
\usepackage{amsfonts}
\usepackage{color,xcolor}
\usepackage[left=2.0cm,right=2.0cm,top=2.0cm,bottom=2.0cm]{geometry}
\usepackage[colorlinks,citecolor=blue,urlcolor=blue]{hyperref}

\newtheorem{theorem}{Theorem}[section]

\newtheorem{lemma}{Lemma}[section]
\newtheorem{proposition}{Proposition}[section]

\newtheorem{definition}{Definition}[section]

\def\norm#1{\|#1\|}
\newcommand{\bal}{\begin{align}}
\newcommand{\bbal}{\begin{align*}}
\newcommand{\beq}{\begin{equation}}
\newcommand{\eeq}{\end{equation}}
\newcommand{\bca}{\begin{cases}}
\newcommand{\eca}{\end{cases}}
\def\div{\mathord{{\rm div}}}
\newcommand{\pa}{\partial}
\newcommand{\fr}{\frac}
\newcommand{\na}{\nabla}
\newcommand{\De}{\Delta}

\newcommand{\cd}{\cdot}

\newcommand{\dd}{\mathrm{d}}

\newcommand{\R}{\mathbb{R}}

\newcommand{\Z}{\mathbb{Z}}

\newcommand{\les}{\lesssim}

\newcommand{\bi}{\Big}
\newcommand{\g}{\big}
\begin{document}
\title{Non-uniform dependence on initial data for the Euler equations in Besov spaces}

\author{Jinlu Li$^{1}$, Yanghai Yu$^{2}$ and Weipeng Zhu$^{3}$\footnote{E-mail: lijl29@gnnu.edu.cn; yuyanghai214@sina.com(Corresponding author); mathzwp2010@163.com}\\
\small $^1$\it School of Mathematics and Computer Sciences, Gannan Normal University, Ganzhou 341000, China\\
\small $^2$\it School of Mathematics and Statistics, Anhui Normal University, Wuhu, Anhui, 241000, China\\
\small $^3$\it School of Mathematics and Information Science, Guangzhou University, Guangzhou 510006, China}

\date{}

\maketitle\noindent{\hrulefill}

{\bf Abstract:} In the paper, we consider the initial value problem to the higher dimensional Euler equations in the whole space. Based on the new technical which is developed in \cite{Li2}, we proved that the data-to-solution map of this problem is not uniformly continuous in nonhomogeneous Besov spaces in the sense of Hadamard. Our obtained result improves considerably the recent result given by Pastrana \cite{Pastrana}.


\vskip0mm\noindent{\hrulefill}

\section{Introduction}

In this article, we consider the Euler equations governing the motion of an incompressible
fluid
\begin{align}\label{eq:E}
(\rm{\mathbf{E}})\quad\begin{cases}
\pa_t u+u\cdot \nabla u+\nabla P=0, &\quad (t,x)\in \R^+\times\R^d,\\
\mathrm{div\,} u=0,&\quad (t,x)\in \R^+\times\R^d,\\
u(0,x)=u_0, &\quad x\in \R^d,
\end{cases}
\end{align}
where the vector field $u(t,x):[0,\infty)\times {\mathbb R}^d\to {\mathbb R}^d$ stands for the velocity of the fluid, the quantity $P(t,x):[0,\infty)\times {\mathbb R}^d\to {\mathbb R}$ denotes
the scalar pressure, and $\mathrm{div\,} u=0$ means that the fluid is incompressible. The mathematical study of the Euler equations of ideal hydrodynamics has a long and distinguished history. We do not detail the literature since it is huge and refer the readers to see
the monographs of Majda--Bertozzi \cite{Majda} and Bahouri--Chemin--Danchin \cite{B.C.D} for
fundamental results and additional references.

The continuous dependence is particularly important when PDEs are used to model
phenomena in the natural world since measurements are always
associated with errors. One of the first results of this type was proved by Kato \cite{Kato} who showed that the
solution operator for the (inviscid) Burgers equation is not H\"{o}lder continuous in the
$H^s(\mathbb{T})$-norm $(s > 3/2)$ for any H\"{o}lder exponent. After the phenomenon of non-uniform continuity for some dispersive equations was studied by Kenig et.al. \cite{Kenig2001}, many results with regard to the non-uniform dependence on the initial data have been obtained, see for example Koch--Tzvetkov \cite{Koch} for the Benjamin-Ono equation, Himonas et.al. \cite{H-K,H-K-M} for the Camassa-Holm equation, Holmes--Keyfitz--Tiglay \cite{Holmes1} for compressible gas in the Sobolev spaces, and Holmes--Tiglay \cite{Holmes2} for the Hunter Saxton equation in Besov spaces, and Himonas--Misio{\l}ek \cite{HM} for the Euler equations.

Next, we mainly recall some of the recent progress which are closely related to our problem. Cheskidov--Shvydkoy
\cite{Cheskidov} proved that the solution to $(\mathbf{E})$ cannot be continuous in the spaces $B^s_{r,\infty}(\mathbb{T}^d)$. Furthermore, Bourgain and Li in \cite{Bourgain1,Bourgain2} proved strong local ill-posedness of $(\mathbf{E})$ in borderline Besov spaces $B^{d/p+1}_{p,r}$ with $(p,r)\in[1,\infty)\times(1,\infty]$ when $d=2,3$. Here we mentioned that the beautiful results of Himonas and Misio{\l}ek \cite{HM} covered both the torus $\mathbb{T}^d$ and the whole spaces $\R^d$  cases. More precisely, they will prove that the solution map for
system $(\mathbf{E})$ in bi(tri)-dimension is not uniformly continuous on bounded
sets into $\mathcal{C}([0, T ], H^s(\mathbb{T}^d\;\rm{or}\; \R^d))$ for any $s\in \R$. Liu--Tang \cite{Liu} extended the periodic result in \cite{HM} to
$B^s_{2,\infty}(\mathbb{T}^d)$ for $r\in[1,\infty]$ and proved that the solution map is not globally uniformly continuous, which is further
extended by Pastrana \cite{Pastrana} to general Besov Spaces $B^s_{p,r}(\mathbb{T}^d)$ for $(p,r)\in[1,\infty]^2$. Meanwhile, for the non-periodic case, Pastrana also obtained the following non-uniform continuous result:
\begin{theorem}[Theorem 2, \cite{Pastrana}]\label{the0}
Let $d=2,3$ and $s>\frac{d}{2}+1, r\in [2,\infty]$. The system $(\mathbf{E})$ is not uniformly continuous from any bounded subset in $B^s_{2,r}$ into $\mathcal{C}([0,T];B^s_{2,r})$.
\end{theorem}
In this paper, we consider the the property of continuous dependence of solutions of the Cauchy problem for system $(\mathbf{E})$. Motivated by our recent work \cite{Li1,Li2}, we generalize the above result to the case $p\neq2$. Now we state main result as follows

\begin{theorem}\label{th1.1}
Let $d\geq 2$. Assume that $(s,p,r)$ satisfies
\begin{align}\label{eq:spr}
s>\frac{d}{p}+1, p\in [1,\infty], r\in [1,\infty) \quad   \mathrm{or}    \quad s=\frac{d}{p}+1, p\in [1,\infty), r=1.
\end{align}
Then system $(\mathbf{E})$ is not uniformly continuous from any bounded subset in $B^s_{p,r}$ into $\mathcal{C}([0,T];B^s_{p,r})$. More precisely, there exists two sequences of solutions $\mathbf{S}_t(f_n+g_n)$ and $\mathbf{S}_t(f_n)$ such that
\bbal
&||f_n||_{B^s_{p,r}}\lesssim 1 \quad\text{and}\quad \lim_{n\rightarrow \infty}||g_n||_{B^s_{p,r}}= 0
\end{align*}
but
\bbal
\liminf_{n\rightarrow \infty}||\mathbf{S}_t(f_n+g_n)-\mathbf{S}_t(f_n)||_{B^s_{p,r}}\gtrsim t,  \quad \forall \;t\in[0,T].
\end{align*}
\end{theorem}

\noindent\textbf{Organization of our paper.}
In Section 2, we list some notations and known results which will be used in the sequel. In Section 3, we present the local well-posedness result and establish some technical Lemmas. In Section 4, we prove our main theorem. Here we give an overview of the strategy:
\begin{itemize}
  \item Choosing a sequence of approximate initial data $f_n$, which can approximate the solution $\mathbf{S}_t(f_n)$;
  \item Considering the initial data $u^n_0=f_n+g_n$ (see Section 3.2 for the constructions of $f_n$ and $g_n$), we shall use a completely new idea. Let us make it more precise: we introduce
  \bbal
\mathbf{w}_n=\mathbf{S}_{t}(u^n_0)-u^n_0-t\mathcal{P}(\mathbf{v}_0^{n})\quad\mbox{ with }\;\mathbf{v}^n_0=-u^n_0\cd\na u^n_0,
\end{align*}
based on the special choice of $f_n$ and $g_n$, we make an important observation that the appearance of $g_n\partial_xf_n$ plays an essential role since it would not small when $n$ is large enough;
  \item The key step is to compute the error $\mathbf{w}_n$ and estimate the $B_{p,r}^s$-norm of this error;
  \item With the approximate solutions $\mathbf{S}_{t}(f_n)$ and $\mathbf{S}_{t}(u^n_0)$ were constructed, combining the precious steps, we can conclude that their distance at the initial time is converging to zero, while at any later time it is bounded below by a positive constant, namely,
\bbal
\lim_{n\rightarrow \infty}||f_n+g_n-f_n||_{B^s_{p,r}}=0
\end{align*}
but
\bbal
\liminf_{n\rightarrow \infty}||\mathbf{S}_t(f_n+g_n)-\mathbf{S}_t(f_n)||_{B^s_{p,r}}\gtrsim t\quad\text{for} \ t \ \text{small enough}.
\end{align*}
 That means the solution map is not uniformly continuous.
\end{itemize}

\section{Littlewood-Paley analysis}
We will use the following notations throughout this paper.
\begin{itemize}
  \item Given a Banach space $X$, we denote its norm by $\|\cdot\|_{X}$.
  \item
\noindent The symbol $A\lesssim B$ means that there is a uniform positive constant $c$ independent of $A$ and $B$ such that $A\leq cB$.
  \item Let us recall that for all $u\in \mathcal{S}'$, the Fourier transform $\mathcal{F}u$, also denoted by $\hat{u}$, is defined by
$$
\mathcal{F}u(\xi)=\hat{u}(\xi)=\int_{\R^d}e^{-ix\cd \xi}u(x)\dd x \quad\text{for any}\; \xi\in\R^d.
$$
  \item The inverse Fourier transform allows us to recover $u$ from $\hat{u}$:
$$
u(x)=\mathcal{F}^{-1}\hat{u}(x)=\frac{1}{2\pi}\int_{\R}e^{ix\xi}\hat{u}(\xi)\dd\xi.
$$
 \item Due to the Holdge decomposition, we know that any vector field $f= (f_1,...,f_d)$ with components in $S_h'(\mathbb{R}^d)$ may be decomposed into one potential part $\mathcal{Q}f$ and one divergence-free part $\mathcal{P}f$, where the projectors $\mathcal{P}$ and $\mathcal{Q}$ are defined by
\begin{eqnarray}\label{Equ2.1}
\mathcal{P} := \mathbb{I}+(-\Delta)^{-1}\nabla {\rm{div}} \quad\mbox{and}\quad \mathcal{Q} := -(-\Delta)^{-1}\nabla {\rm{div}}.
\end{eqnarray}
\end{itemize}
Next, we will recall some facts about the Littlewood-Paley decomposition, the nonhomogeneous Besov spaces and their some useful properties (see \cite{B.C.D} for more details).

There exists a couple of smooth functions $(\chi,\varphi)$ valued in $[0,1]$, such that $\chi$ is supported in the ball $\mathcal{B}\triangleq \{\xi\in\mathbb{R}:|\xi|\leq \frac 4 3\}$, and $\varphi$ is supported in the ring $\mathcal{C}\triangleq \{\xi\in\mathbb{R}:\frac 3 4\leq|\xi|\leq \frac 8 3\}$. Moreover,
\begin{eqnarray*}
\chi(\xi)+\sum_{j\geq0}\varphi(2^{-j}\xi)=1 \quad \mbox{ for any } \xi\in \R.
\end{eqnarray*}
It is easy to show that $\varphi\equiv 1$ for $\frac43\leq |\xi|\leq \frac32$.

For every $u\in \mathcal{S'}(\mathbb{R})$, the inhomogeneous dyadic blocks ${\Delta}_j$ are defined as follows
\begin{numcases}{\Delta_ju=}
0, & if $j\leq-2$;\nonumber\\
\chi(D)u=\mathcal{F}^{-1}(\chi \mathcal{F}u), & if $j=-1$;\nonumber\\
\varphi(2^{-j}D)u=\mathcal{F}^{-1}\g(\varphi(2^{-j}\cdot)\mathcal{F}u\g), & if $j\geq0$.\nonumber
\end{numcases}
In the inhomogeneous case, the following Littlewood-Paley decomposition makes sense
$$
u=\sum_{j\geq-1}{\Delta}_ju\quad \text{for any}\;u\in \mathcal{S'}(\mathbb{R}).
$$

\begin{definition}[\cite{B.C.D}]
Let $s\in\mathbb{R}$ and $(p,r)\in[1, \infty]^2$. The nonhomogeneous Besov space $B^s_{p,r}(\R)$ consists of all tempered distribution $u$ such that
\begin{align*}
||u||_{B^s_{p,r}(\R)}\triangleq \Big|\Big|\g(2^{js}||\Delta_j{u}||_{L^p(\R)}\g)_{j\in \Z}\Big|\Big|_{\ell^r(\Z)}<\infty.
\end{align*}
\end{definition}

Next we recall some nonlinear estimates which will be used for the estimate of pressure term.

\begin{lemma}[\cite{Li1}]\label{lem:P}
Assume $(s,p,r)$ satisfies \eqref{eq:spr}. Then

1) there exists a constant $C=C(d,p,r,s)$ such that for all $u,f\in B^s_{p,r}$ with $\mathrm{div\,} u=0$,
\[\|u\cdot \na f\|_{B^{s-1}_{p,r}}\leq C\|u\|_{B^{s-1}_{p,r}}\|f\|_{B^s_{p,r}}.\]

2) there exists a constant $C=C(d,p,r,s)$ such that for all $u,v\in B^s_{p,r}$ with $\mathrm{div\,} u=\mathrm{div\,} v=0$,
\begin{align*}
\|\na(-\De)^{-1}\mathrm{div\,}(u\cdot \na v)\|_{B^s_{p,r}}&\leq C \big(\|u\|_{C^{0,1}}\|v\|_{B^s_{p,r}}+\|v\|_{C^{0,1}}\|u\|_{B^s_{p,r}}\big);\\
\|\na(-\De)^{-1}\mathrm{div\,}(u\cdot \na v)\|_{B^{s-1}_{p,r}}&\leq C \min\big(\|u\|_{B^{s-1}_{p,r}}\|v\|_{B^s_{p,r}},\|v\|_{B^{s-1}_{p,r}}\|u\|_{B^s_{p,r}}\big),
\end{align*}
where $\|f\|_{C^{0,1}}=\|f\|_{L^\infty}+\|\na f\|_{L^\infty}$.
\end{lemma}

We need an estimate for the transport equation. Consider the following equation:
\begin{align}\label{eq:TDep}
\begin{cases}
\partial_t f+v\cdot \nabla f=g,\\
f(x,t=0)=f_0,
\end{cases}
\end{align}
where $v:{\mathbb R}\times {\mathbb R}^d \to {\mathbb R}^d$, $f_0:{\mathbb R}^d\to {\mathbb R}^N$, and $g:{\mathbb R}\times {\mathbb R}^d\to {\mathbb R}^N$ are given.

\begin{lemma}[Theorem 3.38, \cite{B.C.D} and Lemma 2.9, \cite{Li3}]\label{lem:TDe}
Let $1\leq p,r\leq \infty$. Assume that
\begin{align}
\sigma> -d \min(\frac{1}{p}, \frac{1}{p'}) \quad \mathrm{or}\quad \sigma> -1-d \min(\frac{1}{p}, \frac{1}{p'})\quad \mathrm{if} \quad \mathrm{div\,} v=0.
\end{align}
There exists a constant $C=C(d,p,r,\sigma)$ such that for any smooth solution to \eqref{eq:TDep} and $t\geq 0$ we have
\begin{align}\label{ES2}
\sup_{s\in [0,t]}\|f(s)\|_{B^{\sigma}_{p,r}}\leq Ce^{CV_{p}(v,t)}\Big(\|f_0\|_{B^\sigma_{p,r}}
+\int^t_0\|g(\tau)\|_{B^{s}_{p,r}}\dd \tau\Big),
\end{align}
with
\begin{align*}
V_{p}(v,t)=
\begin{cases}
\int_0^t \|\nabla v(s)\|_{B^{\frac{d}{p}}_{p,\infty}\cap L^\infty}\dd s,&\quad\mathrm{if} \; \sigma<1+\frac{d}{p},\\
\int_0^t \|\nabla v(s)\|_{B^{\sigma}_{p,r}}\dd s,&\quad\mathrm{if} \; \sigma=1+\frac{d}{p} \mbox{ and } r>1,\\
\int_0^t \|\nabla v(s)\|_{B^{\sigma-1}_{p,r}}\dd s, &\quad \mathrm{if} \;\sigma>1+\frac{d}{p}\ \mathrm{or}\ \g\{\sigma=1+\frac{d}{p} \mbox{ and } r=1\g\}.
\end{cases}
\end{align*}
If $f=v$, then for all $\sigma>0$ ($\sigma>-1$, if $\mathrm{div\,} v=0$), the estimate \eqref{ES2} holds with
\[V_{p}(t)=\int_0^t \|\nabla v(s)\|_{L^\infty}\dd s.\]
\end{lemma}

\section{Preliminaries}
Before proceeding, we recall the following local well-posedness estimates for the actual solutions.
\subsection{Local well-posedness estimates for the actual solutions}

Let us recall the local well-posedness result for the Euler equations in Besov spaces.

\begin{lemma}[\cite{Li1}]\label{le4}
Assume that $(s,p,r)$ satisfies \eqref{eq:spr} and for any initial data $u_0$ which belongs to $$B_R=\g\{\psi\in B_{p,r}^s: ||\psi||_{B^{s}_{p,r}}\leq R\g\}\quad\text{for any}\;R>0.$$ Then there exists some $T=T(R,s,p,r)>0$ such that the Euler equations has a unique solution $\mathbf{S}_{t}(u_0)\in \mathcal{C}([0,T];B^s_{p,r})$. Moreover, there holds
\begin{align*}
||\mathbf{S}_{t}(u_0)||_{B^s_{p,r}}\leq C||u_0||_{B_{p,r}^s}.
\end{align*}
\end{lemma}

\subsection{Technical Lemmas}
Firstly, we need to introduce smooth, radial cut-off functions to localize the frequency region.

 Let $\hat{\phi}\in \mathcal{C}^\infty_0(\mathbb{R})$ be an even, real-valued and non-negative funtion on $\R$ and satify
\begin{numcases}{\hat{\phi}(x)=}
1, &if $|x|\leq \frac{1}{4^d}$,\nonumber\\
0, &if $|x|\geq \frac{1}{2^d}$.\nonumber
\end{numcases}
Next, we establish the following crucial Lemmas which will be used later on.
\begin{lemma}\label{ley1} For any $p\in[1,\infty]$, then there exists a positive constant $M$ such that
\begin{align}\label{m}
\liminf_{n\rightarrow \infty}\bi\|\phi^2\cos \bi(\fr{17}{12}2^nx\bi)\bi\|_{L^p}\geq M.
\end{align}
\end{lemma}
{\bf Proof}\quad Without loss of generality, we may assume that $p\in[1,\infty)$. By the Fourier iversion formula and the Fubini thereom, we see that
$$||\phi||_{L^\infty}=\sup_{x\in\R}\frac{1}{2\pi}\Big|\int_{\R}\hat{\phi}(\xi)\cos(x\xi)\dd \xi\Big|\leq \frac{1}{2\pi}\int_{\R}\hat{\phi}(\xi)\dd \xi$$
and
$$\phi(0)=\frac{1}{2\pi}\int_{\R}\hat{\phi}(\xi)\dd \xi>0.$$
Since $\phi$ is a real-valued and continuous function on $\R$, then there exists some $\delta>0$ such that
$$\phi(x)\geq \frac{||\phi||_{L^\infty}}{2}\quad\text{ for any }  x\in B_{\delta}(0).$$
Thus, we have
\bbal
\bi\|\phi^2\cos \bi(\fr{17}{12}2^nx\bi)\bi\|^p_{L^p}&\geq \frac{\phi^2(0)}4\int^\delta_{0}\bi|\cos\bi(\fr{17}{12}2^nx\bi)\bi|^p\dd x\\
&=\frac{\delta}{4}\phi^2(0)\frac{1}{2^n\widetilde{\delta}}\int^{2^n\widetilde{\delta}}_{0}|\cos x|^p\dd x\quad\text{with}\;\widetilde{\delta}=\fr{17}{12}\delta.
\end{align*}
Combining the following simple fact
\bbal
\lim_{n\rightarrow \infty}\frac{1}{2^n\widetilde{\delta}}\int_0^{2^n\widetilde{\delta}}|\cos x|^p\dd x=\frac{1}{\pi}\int^\pi_0|\cos x|^p\dd x,
\end{align*}
thus, we obtain the desired result \eqref{m}.
\begin{lemma}\label{ley2} Let $(s,p,r)$ satisfies \eqref{eq:spr}. Define the high frequency function $f_n$ by
\bbal
f_n=2^{-n(s+1)}
\begin{pmatrix}
-\pa_2\\
\pa_1\\
0\\
\cdots\\
0
\end{pmatrix}
\phi(x_1)\cos \bi(\frac{17}{12}2^nx_1\bi)\phi(x_2)\Phi(x_3,\cdots x_n),
\end{align*}
where
\bbal
\Phi(x_3,\cdots x_n)=
\begin{cases}
\phi(x_3)\cdots \phi(x_d), &\quad \mathrm{if} \; d\geq 3,\\
1, &\quad \mathrm{if} \; d=2.
\end{cases}
\end{align*}
Then for any $\sigma\in\R$, we have
\bal\label{y1}
||f_n||_{B^\sigma_{p,r}}\leq C 2^{n(\sigma-s)}||\phi||^d_{L^{p}}.
\end{align}
\end{lemma}
{\bf Proof}\quad Easy computations give that
\bbal
\hat{f}_n=2^{-n(s+1)-1}
\begin{pmatrix}
-i\xi_2\\
i\xi_1\\
0\\
\cdots\\
0
\end{pmatrix}
\bi[\hat{\phi}\bi(\xi_1-\frac{17}{12}2^n\bi)+\hat{\phi}\bi(\xi_1+\frac{17}{12}2^n\bi)\bi]
\hat{\phi}(\xi_2)\hat{\Phi}(\xi_3,\cdots,x_d),
\end{align*}
which implies
\bbal
\mathrm{supp} \ \hat{f}_n\subset \Big\{\xi\in\R: \ \frac{17}{12}2^n-\fr12\leq |\xi|\leq \frac{17}{12}2^n+\fr12\Big\},
\end{align*}
then, we deduce
\begin{numcases}{\Delta_j(f_n)=}
f_n, &if $j=n$,\nonumber\\
0, &if $j\neq n$.\nonumber
\end{numcases}
Thus, the definition of the Besov space tells us that the desired result \eqref{y1}.
\begin{lemma}\label{ley3} Define the low frequency function $g_n$ by
\bbal
g_n=\frac{12}{17}2^{-n}\Phi(x_3,\cdots,x_d)
\begin{pmatrix}
-\phi(x_1)\phi'(x_2)\\
\phi'(x_1)\phi(x_2)\\
0\\
\cdots\\
0
\end{pmatrix}
,\qquad n\gg1.
\end{align*}
Then there exists a positive constant $\widetilde{M}$ such that
\bbal
\liminf_{n\rightarrow \infty}||g_n\cd\na f_n||_{B^s_{p,\infty}}\geq \widetilde{M}.
\end{align*}
\end{lemma}
{\bf Proof}\quad Notice that
\bbal
\mathrm{supp} \ \hat{g}_n\subset \Big\{\xi\in\R: \ 0\leq |\xi|\leq \fr12\Big\}.
\end{align*}
Then, we have
\bbal
\mathrm{supp}\ \widehat{g_n\cd \na f_n}\subset \Big\{\xi\in\R: \ \frac{17}{12}2^n-1\leq |\xi|\leq \frac{17}{12}2^n+1\Big\},
\end{align*}
which implies
\begin{numcases}
{\Delta_j\g(g_n\cd\na f_n\g)=}
g_n\cd \na f_n, &if $j=n$,\nonumber\\
0, &if $j\neq n$.\nonumber
\end{numcases}
By the definitions of $f_n$ and $g_n$, we obtain $\g(g_n\cd \na f_n\g)_i=0$ with $i=3,\cdots, d$. Notice that
\bbal
\g(g_n\cd \na f_n\g)_1&=\frac{12}{17}2^{-n(s+2)}\Phi^2(x_3,\cdots x_d)\phi(x_1)\phi'(x_1)\cos \bi(\frac{17}{12}2^nx_1\bi)[\phi'(x_2)]^2
\\&\quad -2^{-n(s+1)}\Phi^2(x_3,\cdots x_d)[\phi(x_1)]^2\sin \bi(\frac{17}{12}2^nx_1\bi)[\phi'(x_2)]^2
\\&\quad -\frac{12}{17}2^{-n(s+2)}\Phi^2(x_3,\cdots x_d)\phi'(x_1)\phi(x_2)\phi''(x_2)\phi(x_1)\cos \bi(\frac{17}{12}2^nx_1\bi),
\end{align*}
and
\bbal
\g(g_n\cd \na f_n\g)_2&=\frac{17}{12}2^{-ns}\Phi^2(x_3,\cdots x_d)\phi^2(x_1)\cos \bi(\frac{17}{12}2^nx_1\bi)\phi(x_2)\phi'(x_2)
\\&\quad +2^{-n(s+1)}\Phi^2(x_3,\cdots x_d)\phi(x_1)\phi'(x_1)\sin \bi(\frac{17}{12}2^nx_1\bi)\phi(x_2)\phi'(x_2)
\\&\quad -\frac{12}{17}2^{-n(s+2)}\Phi^2(x_3,\cdots x_d)\phi(x_1)\phi''(x_1)\cos \bi(\frac{17}{12}2^nx_1\bi)\phi(x_2)\phi'(x_2)
\\&\quad +\frac{12}{17}2^{-n(s+2)}\Phi^2(x_3,\cdots x_d)[\phi'(x_1)]^2\cos \bi(\frac{17}{12}2^nx_1\bi)\phi(x_2)\phi'(x_2).
\end{align*}
Then, we have
\bbal
||\g(g_n\cd \na f_n\g)_1||_{B^s_{p,\infty}}&=2^{ns}||\De_{n}\g(\g(g_n\cd \na f_n\g)_1)||_{L^p}=2^{ns}||\g(g_n\cd \na f_n\g)_1||_{L^p}\leq C2^{-n},
\end{align*}
and
\bbal
||\g(g_n\cd \na f_n\g)_2||_{B^s_{p,\infty}}&=2^{ns}||\De_{n}\g(\g(g_n\cd \na f_n\g)_2)||_{L^p}=2^{ns}||\g(g_n\cd \na f_n\g)_2||_{L^p}
\\&\geq \frac{17}{12}\bar{M}\bi\|\phi^2(x)\cos \bi(\frac{17}{12}2^nx\bi)\bi\|_{L^p(\R)}\bi\|\phi(x)\phi'(x)\bi\|_{L^p(\R)}-C2^{-n},
\end{align*}
where
\begin{numcases}
{\bar{M}=}
1, &if $d=2$,\nonumber\\
||\phi||^{2(d-2)}_{L^{2p}(\R)}, &if $d\geq 3$.\nonumber
\end{numcases}
Thus, the result of Lemma \ref{ley1} enables us to finish the proof of Lemma \ref{ley3}.
\section{Non-uniform continuous dependence}
In this section, we will give the proof of Theorem \ref{th1.1}.
Firstly, based on the special choice of $f_n$, we construct approximate solutions $\mathbf{S}_{t}(f_n)$ to Euler equations, then estimate the error of approximate solutions $\mathbf{S}_{t}(f_n)$ and the initial data $f_n$.
\begin{proposition}\label{pro1}
Under the assumptions of Theorem \ref{th1.1}, then we have for $k=\pm1$
\begin{equation}\label{11}
||\mathbf{S}_{t}(f_n)||_{B^{s+k}_{p,r}}\leq C2^{kn},
\end{equation}
and
\bal\label{l2}
||\mathbf{S}_{t}(f_n)-f_n||_{B^{s}_{p,r}}\leq C2^{-\frac12n(s-1)}.
\end{align}
\end{proposition}
{\bf Proof}\quad The local well-posedness result (see Lemma \ref{le4}) tells us that the approximate solution $\mathbf{S}_{t}(f_n)\in \mathcal{C}([0,T];B^s_{p,r})$ and has common lifespan $T\thickapprox1$. Moreover, there holds
\bal\label{bound:f}
||\mathbf{S}_{t}(f_n)||_{L^\infty_T(B^s_{p,r})}\leq C.
\end{align}
Notice that
\bbal
\na P=-\mathcal{Q}\big(\mathbf{S}_{t}(f_n)\cd\na \mathbf{S}_{t}(f_n)\big),
\end{align*}
then we have
\bbal
||\na P||_{B^{s-1}_{p,r}}\leq C\norm{\mathbf{S}_{\tau}(f_n)}_{B^s_{p,r}}\norm{\mathbf{S}_{\tau}(f_n)}_{B^{s-1}_{p,r}},
\end{align*}
and
\bbal
||\na P||_{B^{s+1}_{p,r}}\leq C\norm{\mathbf{S}_{\tau}(f_n)}_{B^s_{p,r}}\norm{\mathbf{S}_{\tau}(f_n)}_{B^{s+1}_{p,r}}.
\end{align*}
By Lemmas \ref{lem:P}-\ref{lem:TDe} and \eqref{bound:f}, we have for any $t\in[0,T]$ and for $k=\pm1$
\begin{align*}
\|\mathbf{S}_{t}(f_n)\|_{B^{s+k}_{p,r}}&\leq C e^{CV_p(\mathbf{S}_{t}(f_n),t)} \big(\|f_n\|_{B^{s+k}_{p,r}}+\int^t_0\|\nabla P(\tau)\|_{B^{s+k}_{p,r}}d \tau\big)
\\&\leq C\big(\norm{f_n}_{B^{s+k}_{p,r}}+\int^t_0
\norm{\mathbf{S}_{\tau}(f_n)}_{B^s_{p,r}}\norm{\mathbf{S}_{\tau}(f_n)}_{B^{s+k}_{p,r}}d \tau\big)
\end{align*}
and Gronwall's inequality, we obtain
\begin{align*}
\norm{\mathbf{S}_{t}(f_n)}_{B^{s+k}_{p,r}}\leq Ce^{C\int_0^t \norm{\mathbf{S}_{\tau}(f_n)}_{B^s_{p,r}} d\tau}\norm{f_n}_{B^{s+k}_{p,r}}\leq C\norm{f_n}_{B^{s+k}_{p,r}}.
\end{align*}
which leads to
\bal\label{y2}
||\mathbf{S}_{t}(f_n)||_{B^{s-1}_{p,r}}\leq C2^{-n} \quad\text{and}\quad ||\mathbf{S}_{t}(f_n)||_{B^{s+1}_{p,r}}\leq C2^{n}.
\end{align}
Setting $\mathbf{\widetilde{u}}=\mathbf{S}_{t}(f_n)-f_n$, then we deduce from \eqref{eq:E} that
\bbal
\pa_t\mathbf{\widetilde{u}}+\mathbf{S}_{t}(f_n)\cd\na\mathbf{\widetilde{u}}=-\mathbf{\widetilde{u}}\cd\na f_n-f_n\cd\na f_n-\na P',\quad \mathbf{\widetilde{u}}_0=0.
\end{align*}
Utilizing Lemma \ref{lem:TDe} and \eqref{bound:f} yields
\bal\label{4.11}
||\mathbf{\widetilde{u}}||_{B^{s-1}_{p,r}}&\leq C e^{CV_p(\mathbf{S}_{t}(f_n),t)}  \bi(\int^t_0\g\|\mathbf{\widetilde{u}}\cd\na f_n,\na P'\g\|_{B^{s-1}_{p,r}}\dd \tau
+t||f_n\cd\na f_n||_{B^{s-1}_{p,r}}\bi)\nonumber
\\&\leq C\bi(\int^t_0\g\|\mathbf{\widetilde{u}}\cd\na f_n,\na P'\g\|_{B^{s-1}_{p,r}}\dd \tau
+t||f_n\cd\na f_n||_{B^{s-1}_{p,r}}\bi),
\end{align}
Since
$$-\na P'=\mathcal{Q}\g(\mathbf{S}_{t}(f_n)\cd\na\mathbf{\widetilde{u}}+\mathbf{\widetilde{u}}\cd\na f_n+f_n\cd\na f_n\g),$$
then by Lemma \ref{lem:P}
$$||\na P'||_{B^{s-1}_{p,r}}\leq ||f_n\cd\na f_n||_{B^{s-1}_{p,r}}+||\mathbf{\widetilde{u}}||_{B^{s-1}_{p,r}}
||f_n,\mathbf{S}_{t}(f_n)||_{B^s_{p,r}}.$$
Combining  Lemma \ref{lem:P} yields
\bbal
&||\mathbf{\widetilde{u}}\cd\na f_n||_{B^{s-1}_{p,r}}\leq C||\mathbf{\widetilde{u}}||_{B^{s-1}_{p,r}}||f_n||_{B^s_{p,r}},\\
&||f_n\cd\na f_n||_{B^{s-1}_{p,r}}\leq ||f_n||_{L^\infty}||f_n||_{B^s_{p,r}}+||\pa_xf_n||_{L^\infty}||f_n||_{B^{s-1}_{p,r}}\leq C2^{-sn}.
\end{align*}
Plugging the above inequalities into \eqref{4.11}, then by \eqref{y1} and \eqref{bound:f}, we infer
\bbal
||\mathbf{\widetilde{u}}||_{B^{s-1}_{p,r}}\leq C   \bi(\int^t_0||\mathbf{\widetilde{u}}||_{B^{s-1}_{p,r}}\dd \tau
+2^{-sn}\bi),
\end{align*}
which along with Gronwall's inequality implies
\bal\label{4.12}
||\mathbf{S}_{t}(f_n)-f_n||_{B^{s-1}_{p,r}}\leq C2^{-sn}.
\end{align}
Applying the interpolation inequality, we obtain from \eqref{y1}, \eqref{y2} and \eqref{4.12}
\bbal
||\mathbf{S}_{t}(f_n)-f_n||_{B^{s}_{p,r}}&\leq ||\mathbf{S}_{t}(f_n)-f_n||^\frac12_{B^{s-1}_{p,r}}
||\mathbf{S}_{t}(f_n)-f_n||^\frac12_{B^{s+1}_{p,r}}\leq C2^{-\frac12n(s-1)}.
\end{align*}
Thus we have finished the proof of the Proposition \ref{pro1}.

To obtain the non-uniformly continuous dependence property for the Euler equations, we need to construct a sequence of initial data $u^n_0=f_n+g_n$, which can not approximate the solution $\mathbf{S}_T(u^n_0)$.
\begin{proposition}\label{pro2}
Under the assumptions of Theorem \ref{th1.1}, then we have
\bal\label{l4}
||\mathbf{S}_{t}(u^n_0)-u^n_0-t\mathcal{P}(\mathbf{v}_0^{n})||_{B^{s}_{p,r}}\leq Ct^{2},
\end{align}
where we denote $\mathbf{v}^n_0=-u^n_0\cd\na u^n_0.$
\end{proposition}
{\bf Proof}\quad Obviously, we obtain from Lemmas \ref{ley2}--\ref{ley3} that
\bbal
||u^n_0||_{B^{s+k}_{p,r}}\leq C2^{kn} \quad \text{for }\; k\in\{0,\pm1,2\}.
\end{align*}
Then, Proposition \ref{pro1} directly tells us that for $k=\pm1$
\bal\label{bound:u}
||\mathbf{S}_{t}(u^n_0)||_{B^{s+k}_{p,r}}\leq C2^{kn}.
\end{align}
Next, we can rewrite the solution $\mathbf{S}_{t}(u^n_0)$ as follows:
\bbal
\mathbf{S}_{t}(u^n_0)=u^n_0+t\mathcal{P}(\mathbf{v}_0^{n})+\mathbf{w}_n\quad\mbox{ with }\;\mathbf{v}^n_0=-u^n_0\cd\na u^n_0.
\end{align*}
Using Lemma \ref{lem:P} and the fact that $B^{s-1}_{p,r}(\R)\hookrightarrow L^\infty(\R)$, we have
\bal\label{re3}
||\mathbf{v}^n_0||_{B^{s-1}_{p,r}}&\leq  \nonumber C||u^n_0||_{B^{s-1}_{p,r}}||u^n_0||_{B^{s}_{p,r}}\leq C2^{-n},\\
||\mathbf{v}^n_0||_{B^{s+1}_{p,r}}&\leq ||u^n_0||_{L^\infty}||u^n_0||_{B^{s+2}_{p,r}}+||\pa_xu^n_0||_{L^\infty}||u^n_0||_{B^{s+1}_{p,r}}\\
&\leq C2^{-n}2^{2n}+C2^n\leq C2^n.\nonumber
\end{align}
Similarly, we have
\bal\label{re3-1}
||\mathcal{Q}(\mathbf{v}^n_0)||_{B^{s-1}_{p,r}}\leq C2^{-n},\qquad 
||\mathcal{Q}(\mathbf{v}^n_0)||_{B^{s+1}_{p,r}}\leq C2^n.
\end{align}
Note that $\mathbf{w}_n=\mathbf{S}_{t}(u^n_0)-u^n_0-t\mathcal{P}(\mathbf{v}_0^{n})$, then we can deduce that $\mathbf{w}_n$ satisfy the following equation
\begin{eqnarray}\label{er}
\left\{\begin{array}{ll}
\pa_t\mathbf{w}_n+\mathbf{S}_{t}(u^n_0)\cd\na\mathbf{w}_n=-t\g(u^n_0\cd\na\mathcal{P}(\mathbf{v}_0^{n})
+\mathcal{P}(\mathbf{v}_0^{n})\cd\na u^n_0\g)
-t^2\g(\mathcal{P}(\mathbf{v}_0^{n})\cd\na\mathcal{P}(\mathbf{v}_0^{n})\g)\\
~~~~~~~~~~~~~~~~~~~~~~~~~~~~~\quad-\mathbf{w}_n\cd\na(u^n_0+t\mathcal{P}(\mathbf{v}_0^{n}))-\na P'',\\
\div\ \mathbf{w}_n=0,\\
\mathbf{w}_n(t=0,x)=0,\end{array}\right.
\end{eqnarray}
where
\bbal
-\na P''&=\mathcal{Q}\big(\mathbf{S}_{t}(u^n_0)\cd\na\mathbf{w}_n
+tu^n_0\cd\na\mathcal{P}(\mathbf{v}_0^{n})+t\mathcal{P}(\mathbf{v}_0^{n})\cd\na u^n_0 +t^2\mathcal{P}(\mathbf{v}_0^{n})\cd\na\mathcal{P}(\mathbf{v}_0^{n})
\\&\quad +\mathbf{w}_n\cd\na
(u^n_0+t\mathcal{P}(\mathbf{v}_0^{n}))\big).
\end{align*}
According to Lemma \ref{lem:P}, we have
\bbal
||\na P''||_{B^{s-1}_{p,r}}&\leq C||\mathbf{w}_n||_{B^{s-1}_{p,r}}
||u^n_{0},\mathbf{v}^n_0,\mathbf{S}_{t}(u^n_0)||_{B^s_{p,r}}+
t||u^n_0||_{B^{s-1}_{p,r}}||\mathcal{P}(\mathbf{v}_0^{n})||_{B^{s}_{p,r}}
\\&\qquad +t^2\|\mathcal{P}(\mathbf{v}_0^{n})||_{B^{s-1}_{p,r}}\|\mathcal{P}(\mathbf{v}_0^{n})||_{B^{s}_{p,r}},
\end{align*}
and
\bbal
||\na P''||_{B^{s}_{p,r}}&\leq C||\mathbf{w}_n||_{B^{s}_{p,r}}
||u^n_{0},\mathbf{v}^n_0,\mathbf{S}_{t}(u^n_0)||_{B^s_{p,r}}
+t||u^n_0||_{B^{s}_{p,r}}||\mathcal{P}(\mathbf{v}_0^{n})||_{B^{s}_{p,r}}
\\&\qquad +t^2\|\mathcal{P}(\mathbf{v}_0^{n})||_{B^{s}_{p,r}}\|\mathcal{P}(\mathbf{v}_0^{n})||_{B^{s}_{p,r}}.
\end{align*}
The local well-posedness result (see Lemma \ref{ley1}) tells us that the approximate solution $\mathbf{S}_{t}(u_0^n)\in \mathcal{C}([0,T];B^s_{p,r})$ and has common lifespan $T\thickapprox1$. Moreover, there holds

Utilizing Lemmas \ref{lem:P}-\ref{lem:TDe} and \eqref{bound:u} to \eqref{er}, we have for all $t\in[0,T]$ and for $k\in\{-1,0\}$,
\bal\label{yy}
&\quad \ ||\mathbf{w}_n(t)||_{B^{s+k}_{p,r}}\nonumber
\\&\leq Ce^{CV_p(\mathbf{S}_{t}(u^n_0),t)}\bi(\int^t_0||\mathbf{w}_n\cd\na(u^n_0+\tau\mathcal{P}(\mathbf{v}_0^{n})),\na P''||_{B^{s+k}_{p,r}}\dd \tau\nonumber
\\&\quad +Ct^2||u^n_0\cd\na\mathcal{P}(\mathbf{v}_0^{n}),\mathcal{P}(\mathbf{v}_0^{n})\cd\na u_{0}^n||_{B^{s+k}_{p,r}}+Ct^3||\mathcal{P}(\mathbf{v}_0^{n})\cd\na\mathcal{P}(\mathbf{v}_0^{n})||_{B^{s+k}_{p,r}}\bi)\nonumber
\\&\leq  C\int^t_0||\mathbf{w}_n||_{B^{s+k}_{p,r}}||u^n_{0},\mathcal{P}(\mathbf{v}_0^{n}),\mathbf{S}_{\tau}(u^n_0)||_{B^s_{p,r}}\dd \tau+C(k+1)\int^t_0||\mathbf{w}_n||_{B^{s-1}_{p,r}}||u^n_{0},\mathcal{P}(\mathbf{v}_0^{n})||_{B^{s+1}_{p,r}}\dd \tau
\nonumber\\&\quad+Ct^2||u^n_0\cd\na\mathcal{P}(\mathbf{v}_0^{n}),\mathcal{P}(\mathbf{v}_0^{n})\cd\na u_{0}^n||_{B^{s+k}_{p,r}}+Ct^2||u^n_0||_{B^{s+k}_{p,r}}||\mathcal{P}(\mathbf{v}_0^{n})||_{B^{s+k}_{p,r}}
\nonumber\\&\quad+Ct^3||\mathcal{P}(\mathbf{v}_0^{n})\cd\na\mathcal{P}(\mathbf{v}_0^{n})||_{B^{s+k}_{p,r}}
+Ct^3||\mathcal{P}(\mathbf{v}_0^{n})||_{B^{s+k}_{p,r}}||\mathcal{P}(\mathbf{v}_0^{n})||_{B^{s}_{p,r}}.
\end{align}
Next, we need to estimate the above terms one by one.

{\bf \underline{Case $k=-1$}.}\quad From Lemma \ref{lem:P} and \eqref{re3}-\eqref{re3-1}, we have
\bbal
&||u^n_0\cd\na\mathcal{P}(\mathbf{v}_0^{n})||_{B^{s-1}_{p,r}}\leq
C||u^n_0||_{B^{s-1}_{p,r}}||\mathcal{P}(\mathbf{v}_0^{n})||_{B^{s}_{p,r}}\leq C2^{-n},
\\&||\mathcal{P}(\mathbf{v}_0^{n})\cd\na u^n_0||_{B^{s-1}_{p,r}}\leq
C||\mathcal{P}(\mathbf{v}_0^{n})||_{B^{s-1}_{p,r}}||u^n_0||_{B^{s}_{p,r}}\leq C2^{-n},
\\&||\mathcal{P}(\mathbf{v}_0^{n})\cd\na\mathcal{P}(\mathbf{v}_0^{n})||_{B^{s-1}_{p,r}}\leq
C||\mathcal{P}(\mathbf{v}_0^{n})||_{B^{s-1}_{p,r}}||\mathcal{P}(\mathbf{v}_0^{n})||_{B^{s}_{p,r}}\leq C2^{-n}.
\end{align*}
Gathering all the above estimates together with \eqref{yy} yields
\bbal
||\mathbf{w}_n(t)||_{B^{s-1}_{p,r}}\leq C\int^t_0||\mathbf{w}_n(t)||_{B^{s-1}_{p,r}}\dd \tau+Ct^22^{-n}
\end{align*}
which along with Gronwall's inequality leads to
\bbal
||\mathbf{w}_n(t)||_{B^{s-1}_{p,r}}\leq Ct^22^{-n}.
\end{align*}

{\bf \underline{Case $k=0$}.}\quad From Lemma \ref{lem:P} and \eqref{re3}-\eqref{re3-1}, we have
\bbal
&||u^n_0\cd\na\mathcal{P}(\mathbf{v}_0^{n})||_{B^{s}_{p,r}}\les||u^n_0||_{B^{s-1}_{p,r}}||\mathcal{P}(\mathbf{v}_0^{n})||_{B^{s+1}_{p,r}}+||u^n_0||_{B^{s}_{p,r}}||\mathcal{P}(\mathbf{v}_0^{n})||_{B^{s}_{p,r}}
\les1,
\\&||\mathcal{P}(\mathbf{v}_0^{n})\cd\na u^n_0||_{B^{s}_{p,r}}\les||\mathcal{P}(\mathbf{v}_0^{n})||_{B^{s-1}_{p,r}}||u^n_0||_{B^{s+1}_{p,r}}+||\mathcal{P}(\mathbf{v}_0^{n})||_{B^{s}_{p,r}}||u^n_0||_{B^{s}_{p,r}}
\les1,
\\&||\mathcal{P}(\mathbf{v}_0^{n})\cd\na\mathcal{P}(\mathbf{v}_0^{n})||_{B^{s}_{p,r}}\les||\mathcal{P}(\mathbf{v}_0^{n})||_{B^{s-1}_{p,r}}||\mathcal{P}(\mathbf{v}_0^{n})||_{B^{s+1}_{p,r}}+||\mathcal{P}(\mathbf{v}_0^{n})||_{B^{s}_{p,r}}||\mathcal{P}(\mathbf{v}_0^{n})||_{B^{s}_{p,r}}
\les1.
\end{align*}
Gathering all the above estimates together with \eqref{yy} and using the Gronwall inequality yields
\bbal
||\mathbf{w}_n(t)||_{B^{s}_{p,r}}&\leq Ct^2+C\int^t_02^n||\mathbf{w}_n(\tau)||_{B^{s-1}_{p,r}}\dd \tau\leq Ct^2.
\end{align*}
Thus, we completed the proof of the Proposition \ref{pro2}.

With the Propositions \ref{pro1}--\ref{pro2} at our hand, we can prove our main Theorem.

{\bf Proof of Theorem \ref{th1.1}}\quad
Obviously, we have
\bbal
||u^n_0-f_n||_{B^s_{p,r}}=||g_n||_{B^s_{p,r}}\leq C2^{-n},
\end{align*}
which means that
\bbal
\lim_{n\to\infty}||u^n_0-f_n||_{B^s_{p,r}}=0.
\end{align*}
Furthermore, we deduce that
\bal\label{yyh}
\quad \ ||\mathbf{S}_{t}(u^n_0)-\mathbf{S}_{t}(f_n)||_{B^s_{p,r}}=&~||t\mathcal{P}(\mathbf{v}_0^{n})+g_n+f_n-\mathbf{S}_{t}(f_n)+\mathbf{w}_n||_{B^s_{p,r}}\nonumber\\
\geq&~ ||t\mathcal{P}(\mathbf{v}_0^{n})||_{B^s_{p,r}}-||g_n||_{B^s_{p,r}}-||f_n-\mathbf{S}_{t}(f_n)||_{B^s_{p,r}}-||\mathbf{w}_n||_{B^s_{p,r}}\nonumber\\
\geq&~ t||\mathcal{P}(\mathbf{v}_0^{n})||_{B^s_{p,\infty}}-C2^{-n\min\{s-1,1\}}-Ct^{2}.
\end{align}
Notice that
$$
\mathcal{P}(\mathbf{v}_0^{n})=\mathcal{P}(f_n\cd\na f_n+f_n\cd\na g_n+g_n\cd\na g_n)+g_n\cd\na f_n-\mathcal{Q}(g_n\cd\na f_n),
$$
by simple calculation, we obtain
\bbal
||f_n\cd\na f_n||_{B^s_{p,r}}&\leq ||f_n||_{L^\infty}||f_n||_{B^{s+1}_{p,r}}+||\na f_n||_{L^\infty}||f_n||_{B^{s}_{p,r}}\leq C2^{-n(s-1)},\\
||f_n\cd\na g_n||_{B^s_{p,r}}&\leq ||f_n||_{B^s_{p,r}}||g_n||_{B^{s+1}_{p,r}}\leq C2^{-n},\\
||g_n\cd\na g_n||_{B^s_{p,r}}&\leq ||g_n||_{B^s_{p,r}}||g_n||_{B^{s+1}_{p,r}}\leq C2^{-2n},\\
||\mathcal{Q}(g_n\cd\na f_n)||_{B^s_{p,r}}&=||\mathcal{Q}(f_n\cd\na g_n)||_{B^s_{p,r}}\leq ||f_n||_{B^s_{p,r}}||g_n||_{B^{s+1}_{p,r}}\leq C2^{-n},
\end{align*}
where we use the equality $\div(g_n\cd\na f_n)=\div(f_n\cd\na g_n)$. Hence, it follows from \eqref{yyh} and Lemma \ref{ley3} that
\bbal
\liminf_{n\rightarrow \infty}||\mathbf{S}_t(f_n+g_n)-\mathbf{S}_t(f_n)||_{B^s_{p,r}}\gtrsim t\quad\text{for enough small } t.
\end{align*}
This completes the proof of Theorem \ref{th1.1}.

\vspace*{1em}
\noindent\textbf{Acknowledgements.}  J. Li is supported by the National Natural Science Foundation of China (Grant No.11801090). Y. Yu is supported by the Natural Science Foundation of Anhui Province (No.1908085QA05). W. Zhu is partially supported by the National Natural Science Foundation of China (Grant No.11901092) and Natural Science Foundation of Guangdong Province (No.2017A030310634).
%\vspace*{1em}

\end{document}